\documentclass[11pt,twoside]{article}
\usepackage{latexsym,amsmath,amsthm,amssymb,amsfonts}
\usepackage{mathrsfs}
\usepackage{amscd}
\usepackage{cite}

\usepackage[ansinew]{inputenc}
\numberwithin{equation}{section}

\catcode`@=11

\newskip\plaincentering \plaincentering=0pt plus 1000pt minus 1000pt
\def\@plainlign{\tabskip=0pt\everycr={}}
\def\eqalignno#1{\displ@y \tabskip\plaincentering
  \halign to\displaywidth{\hfil$\@lign\displaystyle{##}$\tabskip\z@skip
    &$\@lign\displaystyle{{}##}$\hfil\tabskip\plaincentering
    &\llap{$\@lign##$}\tabskip\z@skip\crcr
    #1\crcr}}
\def\leqalignno#1{\displ@y \tabskip\plaincentering
  \halign to\displaywidth{\hfil$\@lign\displaystyle{##}$\tabskip\z@skip
    &$\@lign\displaystyle{{}##}$\hfil\tabskip\plaincentering
    &\kern-\displaywidth\rlap{$\@lign##$}\tabskip\displaywidth\crcr
    #1\crcr}}
\def\plainLet@{\relax\iffalse{\fi\let\\=\cr\iffalse}\fi}
\def\plainvspace@{\def\vspace##1{\noalign{\vskip##1}}}

\def\intic@{\mathchoice{\hskip5\p@}{\hskip4\p@}{\hskip4\p@}{\hskip4\p@}}
\def\negintic@
 {\mathchoice{\hskip-5\p@}{\hskip-4\p@}{\hskip-4\p@}{\hskip-4\p@}}
\def\intkern@{\mathchoice{\!\!\!}{\!\!}{\!\!}{\!\!}}
\def\intdots@{\mathchoice{\cdots}{{\cdotp}\mkern1.5mu
    {\cdotp}\mkern1.5mu{\cdotp}}{{\cdotp}\mkern1mu{\cdotp}\mkern1mu
      {\cdotp}}{{\cdotp}\mkern1mu{\cdotp}\mkern1mu{\cdotp}}}
\newcount\intno@
\def\iint{\intno@=\tw@\futurelet\next\ints@}
\def\iiint{\intno@=\thr@@\futurelet\next\ints@}
\def\iiiint{\intno@=4 \futurelet\next\ints@}
\def\idotsint{\intno@=\z@\futurelet\next\ints@}
\def\ints@{\findlimits@\ints@@}
\newif\iflimtoken@
\newif\iflimits@
\def\findlimits@{\limtoken@false\limits@false\ifx\next\limits
 \limtoken@true\limits@true\else\ifx\next\nolimits\limtoken@true\limits@false
    \fi\fi}
\def\multintlimits@{\intop\ifnum\intno@=\z@\intdots@
  \else\intkern@\fi
    \ifnum\intno@>\tw@\intop\intkern@\fi
     \ifnum\intno@>\thr@@\intop\intkern@\fi\intop}
\def\multint@{\int\ifnum\intno@=\z@\intdots@\else\intkern@\fi
   \ifnum\intno@>\tw@\int\intkern@\fi
    \ifnum\intno@>\thr@@\int\intkern@\fi\int}
\def\ints@@{\iflimtoken@\def\ints@@@{\iflimits@
   \negintic@\mathop{\intic@\multintlimits@}\limits\else
    \multint@\nolimits\fi\eat@}\else
     \def\ints@@@{\multint@\nolimits}\fi\ints@@@}
\def\Sb{_\bgroup\vspace@
        \baselineskip=\fontdimen10 \scriptfont\tw@
        \advance\baselineskip by \fontdimen12 \scriptfont\tw@
        \lineskip=\thr@@\fontdimen8 \scriptfont\thr@@
        \lineskiplimit=\thr@@\fontdimen8 \scriptfont\thr@@
        \Let@\vbox\bgroup\halign\bgroup \hfil$\scriptstyle
            {##}$\hfil\cr}
\def\endSb{\crcr\egroup\egroup\egroup}
\def\Sp{^\bgroup\vspace@
        \baselineskip=\fontdimen10 \scriptfont\tw@
        \advance\baselineskip by \fontdimen12 \scriptfont\tw@
        \lineskip=\thr@@\fontdimen8 \scriptfont\thr@@
        \lineskiplimit=\thr@@\fontdimen8 \scriptfont\thr@@
        \Let@\vbox\bgroup\halign\bgroup \hfil$\scriptstyle
            {##}$\hfil\cr}
\def\endSp{\crcr\egroup\egroup\egroup}
\def\Let@{\relax\iffalse{\fi\let\\=\cr\iffalse}\fi}
\def\vspace@{\def\vspace##1{\noalign{\vskip##1 }}}
\def\aligned{\,\vcenter\bgroup\plainvspace@\plainLet@\openup\jot\m@th\ialign
  \bgroup \strut\hfil$\displaystyle{##}$&$\displaystyle{{}##}$\hfil\crcr}
\def\endaligned{\crcr\egroup\egroup}
\def\matrix{\,\vcenter\bgroup\plainLet@\plainvspace@
    \normalbaselines
  \m@th\ialign\bgroup\hfil$##$\hfil&&\quad\hfil$##$\hfil\crcr
    \mathstrut\crcr\noalign{\kern-\baselineskip}}
\def\endmatrix{\crcr\mathstrut\crcr\noalign{\kern-\baselineskip}\egroup
                \egroup\,}
\newtoks\hashtoks@
\hashtoks@={#}
\def\format{\crcr\egroup\iffalse{\fi\ifnum`}=0 \fi\format@}
\def\format@#1\\{\def\preamble@{#1}%
  \def\c{\hfil$\the\hashtoks@$\hfil}%
  \def\r{\hfil$\the\hashtoks@$}%
  \def\l{$\the\hashtoks@$\hfil}%
  \setbox\z@=\hbox{\xdef\Preamble@{\preamble@}}\ifnum`{=0 \fi\iffalse}\fi
   \ialign\bgroup\span\Preamble@\crcr}

\def\cases{\left\{\,\vcenter\bgroup\plainvspace@
     \normalbaselines\openup\jot\m@th
      \plainLet@\ialign\bgroup$\displaystyle{##}$\hfil&\quad$\displaystyle{{}##}$\hfil\crcr
      \mathstrut\crcr\noalign{\kern-\baselineskip}}
\def\endcases{\endmatrix\right.}
\newif\iftagsleft@
\tagsleft@true
\def\TagsOnRight{\global\tagsleft@false}
\def\tag#1$${\iftagsleft@\leqno\else\eqno\fi
 \hbox{\def\pagebreak{\global\postdisplaypenalty-\@M}%
 \def\nopagebreak{\global\postdisplaypenalty\@M}\rm(#1\unskip)}%
  $$\postdisplaypenalty\z@\ignorespaces}
\interdisplaylinepenalty=\@M
\def\plainallowdisplaybreak@{\def\allowdisplaybreak{\noalign{\allowbreak}}}
\def\plaindisplaybreak@{\def\displaybreak{\noalign{\break}}}
\def\align#1\endalign{\def\tag{&}\plainvspace@\plainallowdisplaybreak@\plaindisplaybreak@
  \iftagsleft@\plainlalign@#1\endalign\else
   \plainralign@#1\endalign\fi}
\def\plainralign@#1\endalign{\displ@y\plainLet@\tabskip\plaincentering\halign to\displaywidth
     {\hfil$\displaystyle{##}$\tabskip=\z@&$\displaystyle{{}##}$\hfil
       \tabskip=\plaincentering&\llap{\hbox{\rm(##\unskip)}}\tabskip\z@\crcr
             #1\crcr}}
\def\plainlalign@
 #1\endalign{\displ@y\plainLet@\tabskip\plaincentering\halign to \displaywidth
   {\hfil$\displaystyle{##}$\tabskip=\z@&$\displaystyle{{}##}$\hfil
   \tabskip=\plaincentering&\kern-\displaywidth
        \rlap{\hbox{\rm(##\unskip)}}\tabskip=\displaywidth\crcr
               #1\crcr}}

\def\re@#1{\par\hangindent\parindent\indent\llap{#1\enspace}\ignorespaces}
\def\qfootnote#1{\edef\@sf{\spacefactor\the\spacefactor}{}#1\@sf
      \insert\footins{\let\egroup=}\footnotesize 
      \interlinepenalty100 \let\par=\endgraf
        \leftskip=0pt \rightskip=0pt
        \splittopskip=10pt plus 1pt minus 1pt \floatingpenalty=20000
   \smallskip\re@{#1}\bgroup\strut\aftergroup{\strut\egroup}\let\next}
\catcode`\@=\active
\TagsOnRight

\begin{document}
\title{\bf
Lower order eigenvalues of the poly-Laplacian with any order on
spherical domains \footnote{This research is supported by Project of
Henan Provincial department of Sciences and Technology (No.
092300410143), and NSF of Henan Provincial Education department (No.
2009A110010).} }
\author{Guangyue Huang,\ Bingqing Ma
\footnote{The corresponding author. Email: bqma$@$henannu.edu.cn}\\
{\normalsize Department of Mathematics, Henan Normal University,
Xinxiang 453007, Henan}
\\{\normalsize  People's Republic of China}\\
}
\date{}
\maketitle
\begin{quotation}
\noindent{\bf Abstract.}~ We consider the lower order eigenvalues of
poly-Laplacian with any order on spherical domains. We obtain
universal inequalities for them and show that our results are
optimal.\\
{{\bf Keywords}: eigenvalue; poly-Laplacian.} \\
{{\bf Mathematics Subject Classification}: Primary 35P15; Secondary
53C20.}

\end{quotation}

\section{Introduction}

Let $\Omega$ be a connected bounded domain in an $n$-dimensional
complete Riemannian manifold $M$. In this paper, we consider the
Dirichlet eigenvalue problem of the poly-Laplacian with order $p$:
\begin{equation}\label{Intr1}\left\{\begin{array}{ll}
(-\Delta)^p u=\lambda u &\ \ \ {\rm in}\ \Omega,\\
u={\partial u\over\partial\nu}=\cdots={\partial^{p-1}
u\over\partial\nu^{p-1}}=0 & \ \ \ {\rm on}\
\partial\Omega,
\end{array}\right.\end{equation} where $\Delta$ is the Laplacian in $M$ and
$\nu$ denotes the outward unit normal vector field of
$\partial\Omega$. Let $0<\lambda_1\leq\lambda_2\leq \lambda_3\leq
\cdots\rightarrow+\infty$ denote the successive eigenvalues for
\eqref{Intr1}, where each eigenvalue is repeated according to its
multiplicity.

When $p=1$, the eigenvalue problem \eqref{Intr1} is called a {\em
fixed membrane problem}. For $M=\mathbb{R}^2$ and $p=1$,
Payne-P\'{o}lya-Weinberger \cite{ppw2} proved
$$
\lambda_2+\lambda_3\leq6\lambda_1.
$$
In 1993, for general dimensions $n\geq2$, Ashbaugh and Benguria
\cite{ashbaugh93} proved
$$
\sum_{i=1}^n(\lambda_{i+1}-\lambda_1)\leq 4\lambda_1.
$$
Recently, the inequalities of eigenvalues of the fixed membrane
problem have been generalized to some Riemannian manifolds. For the
related research and improvement in this direction, see
\cite{brands64,cheng08,huang09,sun08} and the references therein. In
particular, Sun-Cheng-Yang \cite{sun08} proved that when $M$ is an
$n$-dimensional unit sphere $\mathbb{S}^n(1)$,
\begin{equation}\label{Intr4}
\sum_{i=1}^n(\lambda_{i+1}-\lambda_1)\leq 4\lambda_1 +n^2.
\end{equation}

When $p=2$, the eigenvalue problem \eqref{Intr1} is called a {\em
clamped plate problem}. For $M=\mathbb{R}^n$ and $p\geq2$,
Cheng-Ichikawa-Mametsuka proved in \cite{chengqingming} that
\begin{equation}\label{Intr3}
\sum_{i=1}^n(\lambda_{i+1}-\lambda_1)\leq4p(2p-1)\lambda_1,
\end{equation}
\begin{equation}\label{Intr4}
\sum_{i=1}^n(\lambda_{i+1}^{\frac{1}{p}}
-\lambda_1^{\frac{1}{p}})^{p-1}\leq(2p)^{p-1}\lambda_1^{\frac{p-1}{
p}}.
\end{equation}
Inequalities \eqref{Intr3} and \eqref{Intr4} include two universal
inequalities of the clamped plate problem announced by Ashbauth in
\cite{ash}. When $M$ is a general complete Riemannian manifold, for
$p=2$, Cheng-Huang-Wei \cite{Chenghuangwei} obtained
\begin{equation}\label{Intr5}\sum_{i=1}^n(\lambda_{i+1}-\lambda_1)^{\frac{1}{2}}\leq
\big(4\lambda_1^{\frac{1}{2}}+n^2H_0^2\big
)^{\frac{1}{2}}\big\{(2n+4)\lambda_1^{\frac{1}{2}}+n^2H_0^2\big\}^{\frac{1}{2}},
\end{equation} where $H_0^2$ is a nonnegative constant which depends only on $M$
and $\Omega$. For $M=\mathbb{S}^n(1)$, we have $H_0^2=1$ such that
\eqref{Intr5} becomes the following inequality:
\begin{equation}\label{Intr6}\sum_{i=1}^n(\lambda_{i+1}-\lambda_1)^{\frac{1}{2}}\leq
\big(4\lambda_1^{\frac{1}{2}}+n^2\big
)^{\frac{1}{2}}\big\{(2n+4)\lambda_1^{\frac{1}{2}}+n^2\big\}^{\frac{1}{2}}.
\end{equation} We remark that when $\Omega=\mathbb{S}^n(1)$, it holds that
$\lambda_1=0$ and $\lambda_2=\cdots=\lambda_{n+1}=n^2$. Therefore,
the inequality \eqref{Intr6} becomes equality. Hence, for
$M=\mathbb{S}^n(1)$, the inequality \eqref{Intr6} is optimal.

In the present article, we consider the eigenvalue problem
\eqref{Intr1} with any $p$ when $M$ is a unit sphere
$\mathbb{S}^n(1)$. we obtain the following result:

\par\bigskip \noindent{\bf Theorem.} {\it
Let $\Omega$ be a bounded domain in an $n$-dimensional unit sphere
$\mathbb{S}^n(1)$. Let $\lambda_i$ be the $i$-th eigenvlaue of the
eigenvalue problem \eqref{Intr1}. Then we have
\begin{equation}\label{Intr7}\aligned \sum_{i=1}^n(\lambda_{i+1}-\lambda_1)^{1\over2}
\leq&\left\{\left(\lambda_1^{1\over
p}+n\right)^p-\lambda_1+4[2^p-(p+1)]\lambda_1^{1\over
p}\left(\lambda_1^{1\over p}+n\right)^{p-2}\right\}^{1\over2}\\
&\times\left\{4\lambda_1^{1\over p}+n^2\right\}^{1\over2}.
\endaligned\end{equation}
}


\par\bigskip \noindent{\bf Remark 1.} {\it
For $p=2$, the inequality \eqref{Intr7} becomes the optimal
inequality \eqref{Intr6}.

}

\par\bigskip \noindent{\bf Remark 2.} {\it
For the unit sphere $S^n(1)$, by taking $\Omega=S^n(1)$, we know
$\lambda_1=0$ and $\lambda_2=\cdots=\lambda_{n+1}=n^p$. Hence, the
inequality \eqref{Intr7} becomes equality. Therefore, our result is
optimal.}

\par\bigskip \noindent{\bf Acknowledgement.} {\it
The authors thank Professor Qing-Ming Cheng for his helpful
discussion and support.}

\section{Proof of Theorem}

Let $u_i$ be the orthonormal eigenfunction corresponding to
eigenvalue $\lambda_i$, that is, $$ \left\{\begin{array}{ll}
(-\Delta)^p u_i=\lambda_i u_i &\ \ \ {\rm in}\ \Omega,\\
u_i={\partial u_i\over\partial\nu}=\cdots={\partial^{p-1}
u_i\over\partial\nu^{p-1}}=0 & \ \ \ {\rm on}\
\partial\Omega,\\
\int\limits\limits_{\Omega}u_iu_j=\delta_{ij}.
\end{array}\right.$$
Let $x^1,x^2,\ldots,x^{n+1}$ be the standard Euclidean coordinate
functions of $\mathbb{R}^{n+1}$, then $$\mathbb{S}^n(1)=\left\{
(x^1,x^2,\ldots, x^{n+1})\in\mathbb{R}^{n+1}\ ;\
\sum_{i=1}^{n+1}(x^i)^2=1\right\}.$$ It is well known that $$\Delta
x^i=-n x^i, \ \ \ i=1,2,\ldots,n+1.$$ Assume that $B$ is an
$(n+1)\times (n+1)$-matrix defined by $B=(b_{ij})$, where
$$b_{ij}=\int\limits_{\Omega}x^{i}u_1u_{j+1}.$$
Using the orthogonalization of Gram and Schmidt, we know that there
exist an upper triangle matrix $R=(r_{ij})$ and an orthogonal matrix
$Q=(q_{ij})$ such that $R=QB$, that is,
$$
r_{ij}=\sum\limits_{k=1}^{n+1}q_{ik}b_{kj}=
\sum\limits_{k=1}^{n+1}\int\limits_{\Omega}q_{ik}x^{k}u_1u_j=0,\ \
2\leq j\leq i\leq n+1.
$$
Defining $h_{i}=\sum\limits_{k=1}^{n+1}q_{ik}x^{k}$, one gets
$$
\int\limits_{\Omega}h_{i}u_1u_j
=\sum\limits_{k=1}^{n+1}\int\limits_{\Omega}q_{ik} x^{k}u_1u_j=0,\ \
2\leq j\leq i\leq n+1.
$$

Setting
$$\varphi_{i}=h_{i}u_1-u_1\int\limits_{\Omega}h_iu_1^2.$$ Then
$$\int\limits_{\Omega}\varphi_iu_j=0,\ \ \ \
{\rm for\ any}\ j\leq i.
$$
It follows from Rayleigh-Ritz inequality that
\begin{equation}\label{Ineq4}\lambda_{i+1}\leq{
\int\limits_\Omega \varphi_i(-\Delta)^p\varphi_i\over
\|\varphi_i\|^2},\end{equation} where $\|f\|^2=\int\limits_\Omega
|f|^2$. By a direct calculation, we derive at
\begin{equation}\label{Ineq5}
\aligned \int\limits_\Omega \varphi_i(-\Delta)^p\varphi_i=&
\int\limits_\Omega \varphi_i(-\Delta)^p(h_iu_1)\\
=&\int\limits_\Omega \varphi_i\{[(-\Delta)^p(h_iu_1)-h_i(-\Delta)^pu_1]+\lambda_1h_iu_1\}\\
=&\lambda_1\|\varphi_i\|^2+\int\limits_\Omega
\varphi_i[(-\Delta)^p(h_iu_1)-h_i(-\Delta)^pu_1]\\
=&\lambda_1\|\varphi_i\|^2+\int\limits_\Omega
h_iu_1[(-\Delta)^p(h_iu_1)-h_i(-\Delta)^pu_1]\\
&-\int\limits_\Omega h_iu_1^2\int\limits_\Omega u_1[(-\Delta)^p(h_iu_1)-h_i(-\Delta)^pu_1]\\
=&\lambda_1\|\varphi_i\|^2+\int\limits_\Omega
h_iu_1[(-\Delta)^p(h_iu_1)-h_i(-\Delta)^pu_1].
\endaligned\end{equation}
Defining
$$\nabla^r=\cases
\Delta^{r/2} & \ \  {\rm when }\  r \ {\rm is \ even},\\
\nabla(\Delta^{(r-1)/2}) &\ \  {\rm when }\ r \ {\rm is \ odd}.
\endcases
$$ Then \eqref{Ineq5} can be written as
\begin{equation}\label{add1}
\int\limits_\Omega
\varphi_i(-\Delta)^p\varphi_i=\lambda_1\|\varphi_i\|^2+\|\nabla^p(h_iu_1)\|^2
-\lambda_1\|h_iu_1\|^2.
\end{equation}
Putting \eqref{add1} into \eqref{Ineq4} yields
\begin{equation}\label{Ineq6}
(\lambda_{i+1}-\lambda_1)\|\varphi_i\|^2\leq \|\nabla^p(h_iu_1)\|^2
-\lambda_1\|h_iu_1\|^2.
\end{equation}

One gets from integration by parts that $$\aligned
\int\limits_\Omega u_1h_i\langle\nabla h_i, \nabla
u_1\rangle=&{1\over4}\int\limits_\Omega \langle\nabla (h_i^2),
\nabla (u_1^2)\rangle
=-{1\over4}\int\limits_\Omega u_1^2\Delta(h_i^2)\\
=&-{1\over2}\int\limits_\Omega u_1^2h_i\Delta
h_i-{1\over2}\int\limits_\Omega u_1^2|\nabla h_i|^2.
\endaligned$$
Hence,
\begin{equation}\label{Ineq7}
\aligned
&\int\limits_\Omega\varphi_i\left(\langle\nabla h_i, \nabla
u_1\rangle+{1\over2}u_1\Delta h_i\right)\\
=&\int\limits_\Omega u_1h_i\langle\nabla h_i, \nabla
u_1\rangle+{1\over2}\int\limits_\Omega u_1^2h_i\Delta
h_i-\int\limits_\Omega h_iu_1^2\left(\int\limits_\Omega
u_1\langle\nabla h_i,
\nabla u_1\rangle+{1\over2}\int\limits_\Omega u_1^2\Delta h_i\right)\\
=&\int\limits_\Omega u_1h_i\langle\nabla h_i, \nabla
u_1\rangle+{1\over2}\int\limits_\Omega
u_1^2h_i\Delta h_i\\
=&-{1\over2}\int\limits_\Omega u_1^2|\nabla h_i|^2\\
=&-{1\over2}\|u_1\nabla h_i\|^2.
\endaligned
\end{equation}
By virtue of \eqref{Ineq6} and \eqref{Ineq7}, it is easy to see
\begin{equation} \label{Ineq8}\aligned
(\lambda_{i+1}-\lambda_1)^{1\over2}\|u_1\nabla
h_i\|^2=&-2(\lambda_{i+1}-\lambda_1)^{1\over2}\int\limits_\Omega\varphi_i\left(\langle\nabla
h_i, \nabla u_1\rangle+{1\over2}u_1\Delta h_i\right)\\
\leq&\delta(\lambda_{i+1}-\lambda_1)\|\varphi_i\|^2
+{1\over\delta}\left\|\langle\nabla
h_i, \nabla u_1\rangle+{1\over2}u_1\Delta h_i\right\|^2\\
\leq&\delta\{\|\nabla^p(h_iu_1)\|^2 -\lambda_1\|h_iu_1\|^2\}
+{1\over\delta}\left\|\langle\nabla h_i, \nabla
u_1\rangle+{1\over2}u_1\Delta h_i\right\|^2,
\endaligned
\end{equation} where $\delta$ is a positive constant.
Summing over $i$ from 1 to $n+1$ for \eqref{Ineq8}, one finds that
 \begin{equation} \label{Ineq9}\aligned
\sum_{i=1}^{n+1}(\lambda_{i+1}-\lambda_1)^{1\over2}\|u_1\nabla
h_i\|^2\leq&\delta\sum_{i=1}^{n+1}\{\|\nabla^p(h_iu_1)\|^2 -\lambda_1\|h_iu_1\|^2\}\\
&+{1\over\delta}\sum_{i=1}^{n+1}\left\|\langle\nabla h_i, \nabla
u_1\rangle+{1\over2}u_1\Delta h_i\right\|^2\\
=&\delta\sum_{i=1}^{n+1}\{\|\nabla^p(x_iu_1)\|^2 -\lambda_1\|x_iu_1\|^2\}\\
&+{1\over\delta}\sum_{i=1}^{n+1}\left\|\langle\nabla x_i, \nabla
u_1\rangle+{1\over2}u_1\Delta x_i\right\|^2\\
=&\delta\sum_{i=1}^{n+1}\int\limits\limits_\Omega u_1x_i
\{(-\Delta)^p(u_1x_i)-x_i(-\Delta)^pu_1\}\\
&+{1\over\delta}\sum_{i=1}^{n+1}\left\|\langle\nabla x_i, \nabla
u_1\rangle+{1\over2}u_1\Delta x_i\right\|^2.
\endaligned
\end{equation} Making use of the same method as proof of Lemma 1 in \cite{chenqian}, it
is easy to prove $$\int\limits_{\Omega}|\nabla u_1|^2\leq
\lambda_1^{1\over p}.$$ Thus,
\begin{equation}\label{Ineq10}\aligned
\sum_{i=1}^{n+1}\left\|\langle\nabla x_i, \nabla
u_1\rangle+{1\over2}u_1\Delta x_i\right\|^2
=&\sum_{i=1}^{n+1}\int\limits_{\Omega}\left(\langle\nabla
x_i, \nabla u_i\rangle+{1\over2}u_1\Delta x_i\right)^2\\
=&\sum_{i=1}^{n+1}\int\limits_{\Omega}\left({1\over 4}u_1^2(\Delta
x_i)^2+\langle\nabla x_i, \nabla u_1\rangle^2+{1\over 2}\Delta
x_i\langle\nabla x_i, \nabla (u_1^2)\rangle\right)\\
=&{n^2\over 4}+\int\limits_{\Omega}|\nabla
u_1|^2\\
\leq&{n^2\over 4}+\lambda_1^{1\over p}.\endaligned\end{equation} It
has been shown in \cite{chengis09} (see Proposition 2.2 of
\cite{chengis09}) that
\begin{equation}\label{prop21}
\aligned \sum_{i=1}^{n+1}&\int\limits\limits_\Omega u_1x_i
\{(-\Delta)^p(u_1x_i)-x_i(-\Delta)^pu_1\}\\
&\leq \left(\lambda_1^{1\over p}+n\right)^p-\lambda_1
+4[2^p-(p+1)]\lambda_1^{1\over p}\left(\lambda_1^{1\over
p}+n\right)^{p-2}.\endaligned
\end{equation}
Inserting \eqref{Ineq10} and \eqref{prop21} into \eqref{Ineq9}, we
infer
\begin{equation}\label{Ineq11}\aligned
\sum_{i=1}^{n+1}(\lambda_{i+1}-\lambda_1)^{1\over2}\|u_1\nabla
h_i\|^2\leq&\delta\left\{\left(\lambda_1^{1\over
p}+n\right)^p-\lambda_1+4[2^p-(p+1)]\lambda_1^{1\over
p}\left(\lambda_1^{1\over
p}+n\right)^{p-2}\right\}\\
&+{1\over\delta}\left\{\lambda_1^{1\over p}+{n^2\over 4}\right\}.
\endaligned\end{equation}
Minimizing the right hand side of \eqref{Ineq11} as a function of
$\delta$ by choosing $$\delta=\left(\lambda_1^{1\over p}+{n^2\over
4}\over \left(\lambda_1^{1\over
p}+n\right)^p-\lambda_1+4[2^p-(p+1)]\lambda_1^{1\over
p}\left(\lambda_1^{1\over p}+n\right)^{p-2}\right)^{1\over2},$$ we
obtain
\begin{equation}\label{Ineq12}\aligned
\sum_{i=1}^{n+1}(\lambda_{i+1}-\lambda_1)^{1\over2}\|u_1\nabla
h_i\|^2\leq&\left\{\left(\lambda_1^{1\over
p}+n\right)^p-\lambda_1+4[2^p-(p+1)]\lambda_1^{1\over
p}\left(\lambda_1^{1\over p}+n\right)^{p-2}\right\}^{1\over2}\\
&\times\left\{4\lambda_1^{1\over p}+n^2\right\}^{1\over2}.
\endaligned\end{equation}

By a transformation of coordinates if necessary, for any point $q$,
one gets
$$
|\nabla h_i|^2\leq1\ \ \ \ {\rm for\ any}\ i.
$$
It follows that
\begin{equation}\label{Ineq13}\aligned
&\sum_{i=1}^{n+1}(\lambda_{i+1}-\lambda_1)^{1\over2}|\nabla
h_i|^2\\
=&\sum_{i=1}^n(\lambda_{i+1}-\lambda_1)^{1\over2}|\nabla
h_i|^2+(\lambda_{n+1}-\lambda_1)^{1\over2}|\nabla
h_{n+1}|^2\\
=&\sum_{i=1}^n(\lambda_{i+1}-\lambda_1)^{1\over2}|\nabla
h_i|^2+(\lambda_{n+1}-\lambda_1)^{1\over2}\left(n-\sum_{i=1}^n|\nabla
h_i|^2\right)\\
=&\sum_{i=1}^n(\lambda_{i+1}-\lambda_1)^{1\over2}|\nabla
h_i|^2+(\lambda_{n+1}-\lambda_1)^{1\over2}\sum_{i=1}^n\left(1-|\nabla
h_i|^2\right)\\
\geq&\sum_{i=1}^n(\lambda_{i+1}-\lambda_1)^{1\over2}|\nabla
h_i|^2+\sum_{i=1}^n(\lambda_{i+1}-\lambda_1)^{1\over2}\left(1-|\nabla
h_i|^2\right)\\
=&\sum_{i=1}^n(\lambda_{i+1}-\lambda_1)^{1\over2}.
\endaligned
\end{equation} From \eqref{Ineq12} and \eqref{Ineq13}, we obtain
$$\aligned \sum_{i=1}^n(\lambda_{i+1}-\lambda_1)^{1\over2}
\leq&\left\{\left(\lambda_1^{1\over
p}+n\right)^p-\lambda_1+4[2^p-(p+1)]\lambda_1^{1\over
p}\left(\lambda_1^{1\over p}+n\right)^{p-2}\right\}^{1\over2}\\
&\times\left\{4\lambda_1^{1\over p}+n^2\right\}^{1\over2},
\endaligned$$ which concludes the proof of Theorem.

\end{document}